\documentclass[a4paper,12pt]{article}
\usepackage{amsmath,amsthm,amssymb}
\usepackage[mathscr]{eucal}
\theoremstyle{plain}
\newtheorem{proposition}{Proposition}
\newtheorem{theorem}{Theorem}
\newtheorem{lemma}{Lemma}
\theoremstyle{definition}

\theoremstyle{remark}
\newtheorem{corollary}{Corollary}

\newcommand{\dom}{\mathop{\rm dom}\nolimits}
\newcommand{\ran}{\mathop{\rm ran}\nolimits}
\newcommand{\rank}{\mathop{\rm rank}\nolimits}
\begin{document}
\author{G.Y. Tsyaputa}
\title{Transformation Semigroups with the Deformed Multiplication}
\date{}
\maketitle Pairwise non-isomorphic semigroups obtained from the
finite inverse symmetric semigroup $\mathcal{IS}_n ,$  finite
symmetric semigroup $\mathcal{T}_n$ and bicyclic semigroup by the
deformed multiplication proposed by Ljapin are classified.


\section{Introduction}
In a famous Ljapin's monograph [1, 393] there is the following
problem. Let $\Omega_1$, $\Omega_2$ be arbitrary nonempty sets,
$S$ be a set of maps from $\Omega_1$ to $\Omega_2$. Fix a map
$\alpha:\Omega_2\rightarrow\Omega_1$ and define the new
multiplication of maps from $S$ via:
$\varphi\circ\psi=\varphi\cdot\alpha\cdot\psi,$ where the symbol
of $\cdot$ denotes a usual composition of the maps (we perform
multiplication from left to right). It is easy to verify that the
defined operation is associative. Ljapin proposed to investigate
the properties of this semigroup with respect to the restrictions
applied to $S$ and $\alpha$. In particular, there appears an
interesting case if $\Omega_1=\Omega_2=\Omega$, $S$ is some
transformation semigroup on $\Omega$, and $\alpha\in S$. Further
this case might be generalized to an arbitrary semigroup $S$: for
a fixed $a$ from $S$ define the operation $*_a$ on $S$ via: $x*_a
y=xay$. The set $S$ with this operation is, obviously, a semigroup
which we denote $(S, *_a)$. The operation $*_a$ is called
\textit{the multiplication deformed by} $a$ (or just \textit{the
deformed multiplication}).

In the present paper we give the full classification of pairwise
non-isomorphic semigroups obtained from the finite inverse
symmetric semigroup $\mathcal{IS}_n$ of all partial injections of
an $n-$element set, finite symmetric semigroup $\mathcal{T}_n$ of
all transformations of an $n-$element set, and bicyclic semigroup
by the deformed multiplication. We follow terminology and notation
from~\cite{Art}.

\section{$\mathcal{IS}_n$ with the deformed multiplication}

For a partial transformation $\alpha\in\mathcal{IS}_n$ denote by
$\dom(\alpha)$ the domain of $\alpha$, and denote by
$\ran(\alpha)$ its image. The value $|\ran(\alpha)|$ is called the
rank of $\alpha$ and is denoted by $\rank(\alpha)$.
\begin{lemma} For an arbitrary $\alpha\in\mathcal{IS}_n$
the number of idempotents in $(\mathcal{IS}_n ,*_{\alpha})$ is
equal to $2^{\rank (\alpha)}$.
\end{lemma}
\begin{proof} Let $\rank (\alpha)=k$. The equality
$\varepsilon=\varepsilon *_{\alpha}\varepsilon=
\varepsilon\alpha\varepsilon$ implies that $\dom
(\varepsilon)\subseteq\ran (\alpha)$ and $\ran
(\varepsilon)\subseteq\dom (\alpha)$, moreover, for any $x\in\dom
(\varepsilon)$ we have $\varepsilon (x)=\alpha^{-1}(x)$. Therefore
an idempotent $\varepsilon$ is completely defined by its domain
$\dom (\varepsilon)$. On the other hand for any subset
$A\subseteq\ran(\alpha)$ the element, $\varepsilon_A $,  such that
 $\dom(\varepsilon_A)=A$ and $\varepsilon_A (x)=\alpha^{-1} (x)$
for all $x\in A$, satisfies the equality $\varepsilon_A
\alpha\varepsilon_A =\varepsilon_A $. In other words,
$\varepsilon_A $ is an idempotent in $(\mathcal{IS}_n
,*_{\alpha})$. Hence there is a one-to-one correspondence between
idempotents in $(\mathcal{IS}_n ,*_{\alpha})$ and the subsets of
$\ran (\alpha)$ which completes the proof.
\end{proof}
\begin{theorem}\label{th1}Semigroups $(\mathcal{IS}_n ,*_{\alpha})$ and $(\mathcal{IS}_n ,*_{\beta})$
are isomorphic if and only if $\rank (\alpha)=\rank (\beta)$.
\end{theorem}
\begin{proof}
Lemma 1 provides the necessity of the condition. Conversely, let
$\rank (\alpha)=\rank (\beta)=k$.  Then there exist permutations
$\tau$ and $\pi$ in $\mathcal{S}_n $ such that $\beta
=\tau\alpha\pi$. Define the map $f: (\mathcal{IS}_n
,*_{\alpha})\rightarrow (\mathcal{IS}_n ,*_{\beta}) $ by
$f(\xi)=\pi^{-1}\xi\tau^{-1}$. Obviously, $f$ is bijective,
moreover, for arbitrary $\xi,\eta\in\mathcal{IS}_n$ we have:
\begin{align*}
f(\xi *_{\alpha} \eta)=\pi^{-1}\xi *_{\alpha}
\eta\tau^{-1}=\pi^{-1}\xi\tau^{-1}\tau\alpha\pi\pi^{-1}
\eta\tau^{-1}=\\
=\pi^{-1}\xi\tau^{-1}\beta\pi^{-1}\eta\tau^{-1}=f(\xi)*_{\beta}f(\eta).
\end{align*}
\end{proof}
\begin{corollary}
There are $(n+1)$ pairwise non-isomorphic semigroups obtained from
$\mathcal{IS}_n$ by the deformed multiplication.
\end{corollary}
\begin{proof}Follows from theorem~\ref{th1} and the fact, that the rank of the
element in $\mathcal{IS}_n$ can be equal to any integer from the
interval $[0,n]$.
\end{proof}

\section{$\mathcal{T}_n$ with the deformed multiplication}
Let $\mathcal{T}_n$ be the symmetric semigroup of all
transformations of the set $N=\{1,2,\ldots ,n\}$. We call
\textit{the type of transformation} $a\in\mathcal{T}_n$ a set
$(\alpha_{1} ,\alpha_{2} ,\ldots ,\alpha_n)$, where $\alpha_k$ is
the number of those elements $y\in N$, whose full inverse image
$a^{-1} (y)$ contains exactly $k$ elements. Obviously, $1\cdot
\alpha_1 +2\cdot \alpha_2 +\dots +n\cdot \alpha_n =n$, and the sum
$\alpha_1 +\alpha_2 +\dots +\alpha_n$ is equal to the cardinality
of the image of $a$.
\begin{theorem} Semigroups $(\mathcal{T}_n ,*_a )$ and $(\mathcal{T}_n ,*_b )$
are isomorphic if and only if transformations $a$ and $b$ have the
same type.
\end{theorem}
\begin{proof}
\textsl{Necessity.} Let $(\mathcal{T}_n ,*_a )$ and
$(\mathcal{T}_n ,*_b )$ be isomorphic semigroups. On semigroup
$(\mathcal{T}_n ,*_a )$ define an equivalence relation $\sim_a$ in
the following way: $x\sim_a y$ if and only if $x*_a u=y*_a u$ for
all $u\in \mathcal{T}_n $. In the same vein define an equivalence
relation $\sim_b $ on $(\mathcal{T}_n ,*_b )$. First we prove that
arbitrary isomorphism $\varphi : (\mathcal{T}_n ,*_a)\rightarrow
(\mathcal{T}_n ,*_b)$ is in accordance with these equivalence
relations, that is, $\varphi (x)\sim_b \varphi(y)$ if and only if
$x\sim_a y$.

Indeed, let $\varphi : (\mathcal{T}_n ,*_a)\rightarrow
(\mathcal{T}_n ,*_b)$ be an isomorphism and let $x\sim_a y$. Then
for all $u\in\mathcal{T}_n $ we have $x*_a u=y*_a u$, and further
$\varphi (x)*_b\varphi(u)=\varphi (y)*_b\varphi(u)$. However
$\varphi(u)$ runs over whole set $\mathcal{T}_n $, so
$\varphi(x)\sim_b \varphi(y)$. The inverse map $\varphi^{-1} :
(\mathcal{T}_n ,*_b)\rightarrow (\mathcal{T}_n ,*_a)$ is also
isomorphic, therefore $\varphi(x)\sim_b \varphi (y)$ implies
$x\sim_a y$. Consequently, $x\sim_a y$ if and only if
$\varphi(x)\sim_b\varphi (y)$.

Therefore any isomorphism between $(\mathcal{T}_n ,*_a)$ and
$(\mathcal{T}_n ,*_b)$ maps equivalence classes of the relation
$\sim_a$ into corresponding equivalence classes of $\sim_b $.
Hence for the relations $\sim_a$ and $\sim_b$ the cardinalities
and the numbers of the equivalence classes must be equal. We show
that by the cardinalities of equivalence classes of the relation
$\sim_a$ we can find the type $(\alpha_1, \alpha_2 ,\ldots
,\alpha_n)$ of $a$ uniquely.
\begin{lemma}
$x\sim_a y$ if and only if $xa=ya $.
\end{lemma}
\begin{proof}
Obviously, the equality $xa=ya$ implies $x\sim_a y $. Now, let
$xa\neq ya $. Then there exists $k$ in $N$ such that $(xa)(k)\neq
(ya)(k)$. Chose an element $u$ from $\mathcal{T}_n $, which has
different meanings in points $(xa)(k)$ and $(ya)(k)$. Then $x*_a
u=xau$ and $y*_a u=yau$ have different meanings in $k$. Hence,
$x*_a u\neq y*_a u$ and $x\nsim_a y$.
\end{proof}

Denote by $\rho_a$ the partition of the set $\{1,2,\dots ,n\}$
induced by $a$ (that is, $x$ and $y$ belong to the same block of
the partition $\rho_a$ if and only if $a(x)=a(y)$). Count the
cardinality of the equivalence class
$\overline{x_{0}}=\{x\,|\,xa=x_0 a\}$ of the relation $\sim_a$ for
a fixed element $x_0 $. First consider the element $y:=x_0
a=\left(
\begin{array}{lcccr}
 1 &2& \dots & n \\
 y_1 & y_2 & \dots & y_n \end{array}
 \right)$. Obviously, every $y_i$ belongs to the image of $a$, $i=1,\dots
 ,n$. Denote by $N_a (a_i)$ the block of the partition $\rho_a $, defined by
 the element $a_i$ from the image of
 $a$. By $n_a (a_i)$ we denote the cardinality of this block. The equality
 $xa=y$ is satisfied if and only if $(xa)(i)=y_i $ for every $i$ or, what is the
same, $x(i)\in N_a (y_i)$. So $x(i)$ can be chosen in $n_a (y_i)$
ways. The images of $x$ in different points are chosen
independently, therefore transformation $x$ can be chosen in
\begin{equation} \label{odyn}
\overset{n}{\underset{i=1}{\prod}}
n_a (y_i)
\end{equation}
ways and we have the cardinality of the class $\overline{x_0}$.

Denote by $m$ the least cardinality of the blocks of the partition
$\rho_a $. The cardinality of the equivalence class of the
relation $\sim_a$ is the least if all multipliers in~(\ref{odyn})
are equal to $m$, and this cardinality is $m^n$. Now we find the
number of equivalence classes $\overline{x_0}$ of the relation
$\sim_a$ which have cardinality $m^n $. To make all multipliers
in~(\ref{odyn}) equal to $m$, there should be $n_a (y_i)=m$ for
all $i$. However $|\{t\,|\,n_a (t)=m\}|=\alpha_m$. So every $y_i$
can be chosen in $\alpha_m$ ways. Since the meanings of $y_i$ for
different $i$ are chosen independently, there are $\alpha_{m}^{n}$
different $y=x_0 a$ such that the corresponding class
$\overline{x_0}$ contains $m^n$ elements. Therefore by the
relation $\sim_a$ we may find the index $m$ and the value
$\alpha_m$ of the first non zero
 component of the type $(\alpha_1 ,\alpha_2 ,\dots ,\alpha_n)$ of
$a$.

Now the components $\alpha_l$ for $l>m$ can be found recursively.
Assume that the components $\alpha_1 ,\alpha_2 ,\dots
,\alpha_{l-1}$ are already known. For the relation $\sim_a$ denote
by $C$ the number of equivalence classes of the cardinality
$l\cdot m^{n-1} $. Then $C$ is equal to the number of sets $(i_1
,i_2 ,\dots , i_n)$ where some of $i_1 ,i_2 ,\dots , i_n$ may
coincide in general, such that:
\begin{equation}
\label{dva} n_a (y_{i_1})n_a (y_{i_2})\dots n_a (y_{i_n})=l\cdot
m^{n-1}
\end{equation}
 Since by assumption $\alpha_1 ,\alpha_2 ,\dots
,\alpha_{l-1}$ are known, we may find the number $A$ of those sets
$(i_1 ,i_2 ,\ldots , i_n)$ for which all multipliers in the left
hand side of~(\ref{dva}) are less than $l$. This number equals
$\overset{n}{\underset{k=1}{\prod}} \alpha_{m_k}$,  where $m\leq
m_k <l$, and $m_1 \cdot m_2 \dots\cdot m_n =l\cdot m^{n-1}$. The
number $B$ of all sets $(i_1 ,i_2 ,\dots , i_n)$ such that one of
the multipliers in the left hand side of~(\ref{dva}) is equal to
$l$ and other $(n-1)$ multipliers equal $m$, is
$n\cdot\alpha_l\cdot\alpha_m ^{n-1}$. Therefore $\alpha_l$ can be
found from the equality $A+B=C$.

Applying the same reasoning to semigroup $(\mathcal{T}_n ,*_b)$ we
can find the type $(\beta_1 ,\beta_2 ,\dots ,\beta_n)$ of $b$ via
the cardinalities of the equivalence classes of the relation
$\sim_b $. For isomorphic semigroups $(\mathcal{T}_n ,*_a)$ and
$(\mathcal{T}_n ,*_b)$ the numbers of the equivalence classes of
the same cardinality coincide, and the values $ \alpha_k$,
$\beta_k$ , $k=1,\dots,n$ are defined by the numbers of
equivalence classes. So for all $k$ we have $\alpha_k =\beta_k$,
that is, elements $a$ and $b$ have the same type.

\textsl{Sufficiency.} Let elements $a$ and $b$ have the type
$(\alpha_1 ,\alpha_2 ,\dots ,\alpha_n)$. Then there exist
permutations $\pi$ and $\tau$ in $\mathcal{S}_n $, such that
$b=\tau a\pi $. The map $f:(\mathcal{T}_n ,*_a)\rightarrow
(\mathcal{T}_n ,*_b)$ such that $f(x)=\pi^{-1}x\tau^{-1}$ defines
the isomorphism between $(\mathcal{T}_n ,*_a)$ and $(\mathcal{T}_n
,*_b)$. Indeed, $f$ is bijective and
\begin{align*} f(x *_a
y)=\pi^{-1}x *_a y\tau^{-1}=\pi^{-1}x\tau^{-1}\tau a\pi\pi^{-1}y
\tau^{-1}=\\
=\pi^{-1}x\tau^{-1}b\pi^{-1}y\tau^{-1}=f(x)*_b f(y).
\end{align*}
\end{proof}
\begin{corollary} Let $p(n)$ denote the number of ways in
which we can split positive integer $n$ into non ordered sum of
the natural integers. Then there are $p(n)$ pairwise
non-isomorphic semigroups obtained from $\mathcal{T}_n$ by the
deformed multiplication.
\end{corollary}
\begin{proposition} In $\mathcal{T}_n$ there are
\begin{displaymath}\frac{n!\binom{n}{\alpha_1}\binom{n-\alpha_1}{\alpha_2}\ldots
\binom{n-\overset{n-1}{\underset{i=1}{\sum}}\alpha_{i}}{\alpha_{n}}}
{\overset{n}{\underset{i=1}{\prod}}(i!)^{\alpha_{i}}}
\end{displaymath}
transformations of the type $(\alpha_1,\alpha_2,\dots,\alpha_n)$.
\end{proposition}
\begin{proof}
To define transformation $a\in\mathcal{T}_n$ of the type
$(\alpha_1,\alpha_2,\dots,\alpha_n)$ first we chose in $\alpha_1$
elements in $N$ which have $1-$element inverse images. This can be
done in $\binom{n}{\alpha_1}$ ways. Further from elements which
left we chose $\alpha_2$ elements which have $2-$element inverse
images, and so on. So the image of $a$ can be defined in
$\binom{n}{\alpha_1}\binom{n-\alpha_1}{\alpha_2}\ldots
\binom{n-\overset{n-1}{\underset{i=1}{\sum}}\alpha_{i}}{\alpha_{n}}$
ways. Now we write these elements in some order and we write
elements which were chosen at $k$ step exactly $k$ times. Then
every of $n!$ permutations $i_1,\ldots ,i_n$ of numbers $1,2,\dots
,n$ defines transformation $a=\left(
\begin{array}{lcccr}
 i_1 &i_2& \ldots & i_n \\
 a_1 & a_2 & \ldots & a_n \end{array} \right)$ of the type
$(\alpha_1,\alpha_2,\dots,\alpha_n)$. However in this way every
transformation $a$ is counted for several times, since
permutations shifting elements with the same images define the
equal transformation. Hence to find the number of transformations
of type $(\alpha_1,\alpha_2,\dots,\alpha_n)$ we need to divide the
value $n!\binom{n}{\alpha_1}\binom{n-\alpha_1}{\alpha_2}\ldots
\binom{n-\overset{n-1}{\underset{i=1}{\sum}}\alpha_{i}}{\alpha_{n}}$
by the repetition factor
$\overset{n}{\underset{i=1}{\prod}}(i!)^{\alpha_{i}}$ with which
every transformation is received.
\end{proof}
\section{Bicyclic semigroup with the deformed multiplication}
Bicyclic semigroup is a semigroup $\mathcal{B}=\langle
a,b|ab=1\rangle$. It is known~\cite{Cliff} that $\mathcal{B}$ is
the inverse semigroup and each element in $\mathcal{B}$ can be
uniquely written in the canonical form $b^m a^k $, $m,k\geq 0$.
Moreover, $(b^m a^k)^{-1}=b^k a^m $.
\begin{proposition} For every $\alpha\in\mathcal{B}$, $\alpha = b^m a^k
$, $\{b^{k+i}a^{k+i},i\geq 0\}$ is the set of idempotents in the
deformed semigroup $(\mathcal{B}, *_{\alpha})$. Moreover,
idempotents form the infinite decreasing chain with respect to a
natural partial order on the set of idempotents.
\end{proposition}
\begin{proof}Element $\varepsilon_i =b^{k+i}a^{m+i}$  is an idempotent in semigroup
$(\mathcal{B}, *_{\alpha})$. Really, $\varepsilon_i
*_{\alpha}\varepsilon_i
=b^{k+i}a^{m+i}b^{m}a^{k}b^{k+i}a^{m+i}=b^{k+i}a^{i}b^{i}a^{m+i}=b^{k+i}a^{m+i}=\varepsilon_i
$.

Now let $\varepsilon=b^t a^s$ be an idempotent of $(\mathcal{B},
*_{\alpha})$. Assume, that $s<m$. Then $b^t a^s
=\varepsilon=\varepsilon *_{\alpha}\varepsilon
=b^{t}a^{s}b^{m}a^{k}b^{t}a^{s}=b^{t+m-s}a^{k}b^{t}a^{s}$. To make
the powers of $a$ in the canonical form of the left and right hand
sides of this equality equal, we need $k\leq t$. Then
$b^{t+m-s+t-k}a^s =b^t a^s $, and $2t+m-s-k=t$, $m-s=k-t$. However
under assumption, $m-s>0$ and $k-t\leq 0$, so the latter equality
is impossible. Hence, $s\geq m$. Then $\varepsilon
*_{\alpha}\varepsilon =b^t a^{s-m+k}b^t a^s $, and $s-m+k=t$,
$s-m=t-k$. Denote $i=s-m$. Then $\varepsilon=b^{k+i}a^{m+i}$.

Let $\varepsilon_i =b^{k+i}a^{m+i}$ and $\varepsilon_j
=b^{k+j}a^{m+j}$ be two idempotents and without loss of generality
let $i\geq j$. Then $\varepsilon_i *_{\alpha}\varepsilon_j
=b^{k+i}a^{m+i}b^{m}a^{k}b^{k+j}a^{m+j}=b^{k+i}a^{i}b^{j}a^{m+j}=b^{k+i}a^{m+i}=\varepsilon_i
$. Analogously $\varepsilon_j *_{\alpha}\varepsilon_i
=\varepsilon_i $, therefore $\varepsilon_i \leq\varepsilon_j$ if
and only if $i\geq j $. So the set of idempotents is linearly
ordered, and $\varepsilon_0 =b^{k+0}a^{m+0}=b^{k}a^{m}$ is maximal
idempotent in $(\mathcal{B}, *_{\alpha})$.
\end{proof}
\begin{theorem} For different $\alpha$ and $\beta$
semigroups $(\mathcal{B}, *_{\alpha})$ and $(\mathcal{B},
*_{\beta})$ are not isomorphic.
\end{theorem}
\begin{proof} Let $\alpha=b^m a^k $, $\beta =b^u a^v $. Take idempotent $\varepsilon_i
=b^{k+i}a^{m+i}$, $i\geq 0$ in semigroup $(\mathcal{B},
*_{\alpha}) $ and consider the sets:
\begin{displaymath}
P^{\alpha}_{i}=\{\xi\in\mathcal{B}\,|\,\varepsilon_i
*_{\alpha}\xi\neq\xi\}\quad\text{ and }\quad
Q^{\alpha}_{i}=\{\xi\in\mathcal{B}\,|\,\xi
*_{\alpha}\varepsilon_i\neq\xi\}.
\end{displaymath}
If element $\xi=b^t a^s$ does not belong to $P^{\alpha}_{i}$ then
$\xi=\varepsilon_i *_{\alpha}\xi $, that is,
$b^{k+i}a^{m+i}b^{m}a^{k}b^{t}a^{s}=b^{t}a^{s}$ and
$b^{k+i}a^{k+i}b^{t}a^{s}=b^{t}a^{s}$.

If $k+i>t $ then $b^{k+i}a^{k+i}b^{t}a^{s}=b^{k+i}a^{s+k+i-t}\neq
b^{t}a^{s}$. Therefore $k+i\leq t $. On the other hand, if
$k+i\leq t$ then
$b^{k+i}a^{k+i}b^{t}a^{s}=b^{k+i}b^{t-(k+i)}a^{s}=b^{t}a^{s}$.
Hence, element $\xi=b^t a^s \notin P^{\alpha}_{i}$ if and only if
$k+i\leq t$.

In the same way we can shown that $\xi=b^t a^s \notin
Q^{\alpha}_{i}$ if and only if $m+i\leq s$.

Thus $\xi=b^t a^s \in P^{\alpha}_{i}\cap Q^{\alpha}_{i} $ if and
only if $t<k+i$ and $s<m+i $. Then the cardinality of an
intersection $P^{\alpha}_{i}\cap Q^{\alpha}_{i}$ equals
$|P^{\alpha}_{i}\cap Q^{\alpha}_{i}|=(k+i)(m+i)$.

It is relatively easy to show that by these cardinalities the
powers of the element $\alpha=b^m a^k$ can be found. Indeed,
$|P^{\alpha}_{1}\cap Q^{\alpha}_{1}|-|P^{\alpha}_{1}\cap
Q^{\alpha}_{0}|=(k+1)(m+1)-(k+1)m=km+k+m+1-km-m=k+1 $, and
$|P^{\alpha}_{1}\cap Q^{\alpha}_{1}|-|P^{\alpha}_{0}\cap
Q^{\alpha}_{1}|=(k+1)(m+1)-k(m+1)=m+1 $.

Applying the same reasoning to $(\mathcal{B}, *_{\beta})$,
$\beta=b^{u}a^v$ we can show that $|P^{\beta}_{1}\cap
Q^{\beta}_{1}|-|P^{\beta}_{1}\cap Q^{\beta}_{0}|=v+1$ and
$|P^{\beta}_{1}\cap Q^{\beta}_{1}|-|P^{\beta}_{0}\cap
Q^{\beta}_{1}|=u+1$.

Hence if semigroups $(\mathcal{B}, *_{\alpha})$ and $(\mathcal{B},
*_{\beta})$ are isomorphic then the cardinalities of corresponding
sets must be equal, that is, $u=m$ and $v=k$, and $\alpha =\beta$.
\end{proof}
\begin{corollary}With respect to isomorphism by the deformed multiplication
we get infinitely many different semigroups from the bicycle
semigroup.
\end{corollary}
\begin{theorem} Semigroups $(\mathcal{B}, *_{\alpha})$ and $(\mathcal{B},
*_{\beta})$ are anti-isomorphic if and only if $\alpha$ and
$\beta$ are inverse.
\end{theorem}
\begin{proof} Let semigroups $(\mathcal{B}, *_{\alpha})$ and $(\mathcal{B},
*_{\beta})$ be such that elements $\alpha$ and $\beta$ are inverse
in $\mathcal{B}$. It is known that the bicycle semigroup is
inverse, and for an element $\alpha =b^m a^k$ there is a unique
inverse element written in a canonical form as $\beta=b^k a^m $.

Consider the map $\varphi :(\mathcal{B}, *_{\alpha})\rightarrow
(\mathcal{B}, *_{\beta})$, $\varphi (b^x a^y)=b^y a^x $. Then
$\varphi $ is a bijection and $\varphi (\alpha)=\varphi (b^m
a^k)=b^k a^m =\beta $. Let $b^x a^y$, $b^t a^s $ be arbitrary
elements in $\mathcal{B}$. If $y\geq t$ then $\varphi (b^x a^y b^t
a^s)=\varphi (b^x a^{y-t+s})=b^{y-t+s}a^x=b^s a^t b^y a^x =
\varphi (b^t a^s)\varphi (b^x a^y)$. The analogous equality is
received if $y<t$. First we prove that $\varphi
(b^{x_1}a^{y_1}\cdots b^{x_n}a^{y_n})=b^{y_n}a^{x_n}\cdots
b^{y_1}a^{x_1}$. In fact,
\begin{align*}
  \varphi (b^{x_1}a^{y_1}b^{x_2}a^{y_2}\cdots
b^{x_n}a^{y_n})=\varphi (b^{x_2}a^{y_2}\cdots
b^{x_n}a^{y_n})\varphi(b^{x_1}a^{y_1})=\\
\varphi
(b^{x_2}a^{y_2}\cdots b^{x_n}a^{y_n})b^{y_1}a^{x_1}=\cdots =
b^{y_n}a^{x_n}\cdots b^{y_1}a^{x_1}.
\end{align*}

 Now for arbitrary $\xi=b^{x}a^{y}$, $\eta =b^{u}a^{v}$ in $\mathcal{B}$ we have
\begin{displaymath} \varphi(\xi *_{\alpha} \eta)=\varphi(b^x a^y b^m a^k b^u
a^v)=b^v a^u b^k a^m b^y a^x=\varphi (\eta) *_{\beta}\varphi
(\xi).
\end{displaymath}
 Hence, $\varphi$ is anti-isomorphism between $(\mathcal{B}, *_{\alpha})$
 and $(\mathcal{B},
*_{\beta})$. However the semigroup anti-isomorphic to the given
one is unique with respect to isomorphism, so if semigroups
$(\mathcal{B}, *_{\alpha})$ and $(\mathcal{B}, *_{\beta})$ are
anti-isomorphic then elements $\alpha$ and $\beta$ are inverse.
\end{proof}

Given to the editorial board on 12.05.2003
\end{document}